\documentclass{amsart}

\usepackage[alphabetic]{amsrefs}

\numberwithin{equation}{subsection}

\swapnumbers
\theoremstyle{plain}
\newtheorem{thm}[subsection]{Theorem}

\newtheorem*{thm*}{Theorem}

\theoremstyle{definition}

\theoremstyle{remark}

\newcommand{\F}{\mathbb{F}}

\newcommand{\Z}{\mathbb{Z}}
\newcommand{\Q}{\mathbb{Q}}

\newcommand{\C}{\mathbb{C}}

\renewcommand{\P}{\mathbb{P}}

\newcommand{\sha}{{\hbox to 10pt{\rlap{\hskip2.8pt\vrule
height6pt\hskip1.6pt\vrule height6pt\hskip1.6pt
\vrule height6pt}\hskip1pt\vrule height0.8pt width 8pt\hskip1pt}}}

\newcommand{\into}{\hookrightarrow}

\newcommand{\nodiv}{\not|}
\def\nodiv{\mathrel{\mathchoice{\not|}{\not|}{\kern-.2em\not\kern.2em|}{\kern-.2em\not\kern.2em|}}}

\DeclareMathOperator{\rk}{Rank}

\DeclareMathOperator{\Hom}{Hom}
\DeclareMathOperator{\aut}{Aut}

\DeclareMathOperator{\lcm}{lcm}
\DeclareMathOperator{\gal}{Gal}

\DeclareMathOperator{\en}{End}

\newcommand{\Qbar}{{\overline{\mathbb{Q}}}}

\newcommand{\kbar}{{\overline{k}}}

\newcommand{\RR}{{\mathcal{R}}}

\newcommand{\End}{\en}
\newcommand{\Gal}{\gal}

\theoremstyle{remark}

\newtheorem{sex}[subsubsection]{Example}
\newtheorem{sexs}[subsubsection]{Examples}

\theoremstyle{plain}
\newtheorem{sthm}[subsubsection]{Theorem}
\newtheorem{slemma}[subsubsection]{Lemma}

\begin{document}

\title[Ranks of Jacobians in towers]{Ranks of Jacobians in towers of function fields}

\author{Douglas Ulmer}
\address{School of Mathematics \\ Georgia Institute of Technology \\
Atlanta, GA, 30332, USA}
\email{ulmer@math.gatech.edu}

\author{Yuri G. Zarhin}
\address{Department of Mathematics \\ Penn
State University \\ University Park, PA 16802, USA}
\email{zarhin@math.psu.edu}


\date{July 6, 2010}

\maketitle

\section{Introduction}
Let $k$ be a field of characteristic zero and $k(t)$ the rational
function field over $k$.  Several authors have considered the behavior
of ranks of elliptic curves in towers of function fields such as
$k(t^{1/n})$.  See for example \cite{Shioda86}, \cite{Ulmer02},
\cite{Silverman04}, \cite{Ellenberg06}, \cite{Ulmer07b},
\cite{Ulmer07c}, \cite{Berger08}, and \cite{UlmerDPCT}.  In
particular, it is known that there are elliptic curves $E$ over $k(t)$
such that the rank of $E(K)$ stays bounded as $K$ varies through the
tower of function fields $K=k(t^{1/p^r})$ with $r\ge1$ or the tower
$K=k(t^{1/d})$ with $d\ge1$.

In this paper we combine the rank formula of \cite{UlmerDPCT} with
strong upper bounds on endomorphisms of Jacobians due to the senior
author (e.g., \cite{ZarhinM}) to give many examples of higher
dimensional, absolutely simple Jacobians over $k(t)$ with bounded rank
in towers $k(t^{1/p^r})$.  In many cases we are able to compute the
rank at every layer of the tower.

Our methods give rise to many examples, but to fix ideas we
state the following typical result.

\begin{thm}\label{thm:main}
Let $g_X$ be an integer $\ge2$ and let $X$ be the smooth, proper curve
of genus $g_X$ over $\Q(t)$ with affine plane model
$$ty^2=x^{2g_X+1}-x+t-1.$$
Let $J$ be the Jacobian of $X$.  Then $J$ is absolutely simple and the
$\overline{\Q(t)}/\Qbar$-trace of $J$ is zero.  Moreover, for every 
prime number $p$ and every integer $r\ge0$, we have
$$\rk J(\Qbar(t^{1/p^r}))=2g_X.$$
\end{thm}

\section{Endomorphism algebras of 
superelliptic Jacobians}\label{s:endos}
\subsection{Notation}\label{ss:endo-intro}
In this section, we collect some results on the endomorphism algebras
of certain superelliptic Jacobians.  Throughout the paper, $\kbar$
will denote an algebraic closure of $k$.

If $X$ and $Y$ are abelian varieties over $k$, then we write
$\End(X)$ and $\Hom(X,Y)$ for the corresponding ring and
group of homomorphisms over $\kbar$.

Let $f$ be a non-constant polynomial with coefficients in $k$ and
without multiple roots.  We write $m$ for the degree of $f$,
$\RR_f$ for the set of roots of $f$ in $\kbar$, and $k(\RR_f)$ for the
splitting field of $f$. As usual
$$\gal(f)=\aut(k(\RR_f)/k).$$

For example, consider $f_m=x^m-x-1$.  It is known
\cite{SerreGalois}*{p.~42} that if $k=\Q$ then $\Gal(f_m)$ is the
symmetric group $S_m$.

\subsection{Superelliptic Jacobians}\label{ss:superelliptic-jacobians}
Given a polynomial $f$ of degree $m$ with distinct roots and a
positive integer $d$, let $C_{f,d}$ be the smooth projective curve
over $k$ with affine model
$$z^d=f(x).$$
Let $J_{f,d}$ be the Jacobian of $C_{f,d}$.  In this paper we will
only consider the case where $d=q=p^r$ with $p$ a prime number.  To
avoid trivialities, we always assume that $m>1$ and $q>1$.

We have an obvious projection $C_{f,q}\to C_{f,q/p}$ and an induced
homomorphism of Jacobians $J_{f,q}\to J_{f,q/p}$.  Since $C_{f,q}\to
C_{f,q/p}$ is totally ramified over the zeroes of $f$, the
homomorphism $J_{f,q}\to J_{f,q/p}$ has connected kernel, which we
call the {\it new part\/} of $J_{f,q}$ and which we denote
$J^{(f,q)}$.  It is easy to see that $J_{f,q}$ is isogenous to
$\prod_{i=1}^r J^{(f,p^i)}$.

A simple
application of the Riemann-Hurwitz formula shows that 
$$\dim(J_{f,q})=(m-1)(q-1)/2-(\gcd(q,m)-1)/2$$
and
$$\dim(J^{(f,q)})=\begin{cases}
(m-1)\phi(q)/2&\text{if $q\nodiv m$}\\
(m-2)\phi(q)/2&\text{if $q|m$.}
\end{cases}$$
Note in particular that $J^{(f,q)}$ is
non-zero except when $m\le2$ and $q\le2$.

We write $\zeta_q$ for a primitive $q$-th root in $\kbar$.  Over
$k(\zeta_q)$, we have a natural action of the $q$-th roots of unity on
$C_{f,q}$, namely $(x,z)\mapsto(x,\zeta_qz)$.  The induced action on
Jacobians gives rise to a homomorphism $\Z[\zeta_q]\to\End(J^{(f,q)})$
which is known to be injective whenever $J^{(f,q)}$ is non-zero.  (See
\cite{ZarhinM}*{4.8}, \cite{ZarhinPisa}*{5.8}, and \cite{Xue}*{2.6};
see also \cite{UlmerDPCT}*{7.8.1} for this result in a slightly more
general context.)  We will be interested in various conditions which
guarantee that this homomorphism is an isomorphism.

\subsection{The case $m=2$}\label{CM}
Suppose that $m=2$, i.e., that $f$ is quadratic.  By a linear change
of variable we may take $f$ to have the form $f(x)=x^2-a$ where $a\in
k$ is non-zero.  If $q>2$ (so that $J^{(f,q)}$ is non-zero), we have
$\dim(J^{(f,q)})=\phi(q)/2$.  The inclusion
$\Z[\zeta_q]\into\End(J^{(f,q)})$ exhibits a commutative integral domain
of rank $2\dim(J^{(f,q)})$ as a subalgebra of $\End(J^{(f,q)})$ and so
over $\kbar$, $J^{(f,q)}$ is an abelian variety of CM type.  By
\cite{Shimura98}*{Prop.~3, p.~36}, $J^{(f,q)}$ is isogenous to a power
$A^e$ of a simple abelian variety $A$ of CM type.  It may happen (for
example if $q$ is even) that $e>1$ in which case $\End(J^{(f,q)})$ is
strictly bigger than $\Z[\zeta_q]$.  In any case, it follows that
$J_{f,q}$ is of CM type over $\kbar$.

\subsection{}\label{ss:3-monomials}
An alternative approach to the case $f(x)=x^2-a$ that works more
generally for $f(x)=x^m-a$ is to note that the curve $z^d=x^m-a$ is
obviously covered (over $\kbar$) by the Fermat curve of degree
$\lcm(m,d)$.  Since the Jacobian of the Fermat curve is an abelian
variety of CM type \cite{Schmidt84}*{VI, 1.2 and 1.5}, the same is
true of $J_{f,d}$ for all $d$.

\subsection{The case $m\ge5$}
For $m>2$ it is again not true in general that
$\Z[\zeta_q]\cong\End(J^{(f,q)})$.  However, one expects that for
``sufficiently general'' $f$, $J^{(f,q)}$ should have no endomorphisms
beyond the obvious ones.  We make this precise in several cases,
starting with $m\ge5$.

\begin{sthm}
\label{bign}
 Suppose that $m \ge 5$ and $\Gal(f)=S_m$ or $A_m$. Then
$$\End(J^{(f,q)})=\Z[\zeta_q].$$
In particular, $J^{(f,q)}$ is an absolutely simple abelian variety that is not
of CM type over $\kbar$.
\end{sthm}

\begin{proof} 
  The endomorphism ring was calculated in earlier works of Zarhin and
  Xue: The case when either $p$ does not divide $m$ or $q$ divides $m$
   is proven in \cite{ZarhinM}*{Th.~4.17 on p.~359}.  The remaining  case when 
   $p$ divides $m$ but $q$ does not divide $m$ is proven in \cite{Xue}.  
   Since the endomorphism ring is a
  domain, $J^{(f,q)}$ is absolutely simple.  Since $m\ge5$, the
  dimension of $J^{(f,q)}$ is strictly greater than $\phi(q)/2$ and so
  $J^{(f,q)}$ is not of CM type.
\end{proof}

\subsection{The cases $m=4$ and $m=3$}
These cases (which are not needed for the basic
Theorem~\ref{thm:main}) are more complicated and require additional
hypotheses.

\begin{sthm}\label{n4} 
  Suppose that $p$ is odd, $m=4$, and $\Gal(f)=S_4$. Let $k^{\prime}$
  be the unique subextension of $k(\RR_f)$ which is quadratic over
  $k$. Assume that $k^{\prime}$ does not lie in $k(\zeta_p)$.
  Then the Galois group of $f(x)$ over $k(\zeta_q)$ is still $S_4$ and
  $$\End(J^{(f,q)})=\Z[\zeta_q].$$
  In particular, $J^{(f,q)}$ is an absolutely simple abelian variety
  that is not of CM type over $\kbar$.
\end{sthm}

\begin{proof}
  In light of \cite{ZarhinMZ2}*{Cor.~1.5 on p.~693 and Cor.~5.3 on
    p.~705}, it suffices to check that the Galois group of $f(x)$ over
  $k(\zeta_q)$ remains $S_4$, i.e., that $k(\RR_f)$ and $k(\zeta_q)$ are
  linearly disjoint over $k$.  Since $p$ is odd,
  $\gal(k(\zeta_q)/k)\cong G_2\times G_{odd}$ where $G_2$ is a cyclic
  group of 2-power order and $G_{odd}$ is a cyclic group of odd order.
  Moreover $k(\zeta_q)^{G_{odd}}\subset k(\zeta_p)$.  Since $S_4$ is
  generated by transpositions, it has no non-trivial quotients of odd
  order.  Thus $k(\RR_f)$ and $k(\zeta_q)$ are linearly disjoint if
  and only if $k'$ and $k(\zeta_p)$ are linearly disjoint (which
  in turn occurs if and only if $k'$ and $k(\zeta_q)^{G_{odd}}$
  are linearly disjoint).
\end{proof}

Ramification conditions give a convenient criterion for linear disjointness:

\begin{slemma}
\label{ramification} Let $L/\Q$ be a finite extension. Suppose that a
prime $p$ is unramified in $L$. Then $L$ and $\Q(\zeta_q)$ are linearly
disjoint over $\Q$.
\end{slemma}

\begin{proof}
This follows immediately from the fact that $\Q(\zeta_q)/\Q$ is
totally ramified over $p$.
\end{proof}

\begin{sthm}
\label{n4bis} Suppose that $k=\Q$, $m=4$ and $\Gal(f)=S_4$.  If $p$ is an odd
prime that is unramified in $\Q(\RR_f)$ then
$$\End(J^{(f,q)})=\Z[\zeta_q].$$
In particular, $J^{(f,q)}$ is an absolutely simple abelian variety that is not
of CM type over $\kbar$.
\end{sthm}

\begin{proof}
It follows from Lemma~\ref{ramification} that $\Q(\RR_f)$ and
$\Q(\zeta_p)$ are linearly disjoint over $\Q$.  Theorem~\ref{n4} then
implies that the endomorphism ring is $\Z[\zeta_q]$.  
Since this is a domain, $J^{(f,q)}$ is absolutely simple and since
$\dim J^{(f,q)}>\phi(q)/2$, it is not of CM type.
\end{proof}

\begin{sexs}
Let $k=\Q$ and $f(x)=f_4(x)=x^4-x-1$. Then the discriminant of $f_4(x)$ is
$-283$ (see \cite{ZarhinMZ2}*{p.~693}).  Since $283$ is a prime, it follows that
$\Q(\RR_{f_4})/\Q$ is unramified outside $283$.
 So, if $p$ is odd and $\ne 283$  then for all $q=p^r$
$$\End(J^{(f_4,q)})=\Z[\zeta_q].$$
(The case of $q=p=3$ is Example 1.6 of \cite{ZarhinMZ2}.)

For another example, again take $k=\Q$ and let $f=x^4-x+2$.   The
discriminant of $f$ is $-27\cdot(-1)^4+256\cdot2^3=2021=43\cdot47$.  Reducing mod
3, one checks that $f$ has no roots in $\F_3$ nor in $\F_9$ and so is
irreducible over $\F_3$ and, {\it a fortiori\/}, irreducible over
$\Q$.  Moreover, the Galois group $\Gal(f)$ is isomorphic to $S_4$.
(Indeed, looking at the Frobenius at $p$ for $p=2,3,5$ shows that
$\Gal(f)$ contains a transposition, a 4-cycle, and a 3-cycle.  Since
the only subgroups of $S_4$ generated by a 4-cycle and a transposition are
all of $S_4$ or a group of order 8, we see that $\Gal(f)$ must be $S_4$.)
The argument of the first example applies when $p$ is
odd and $\neq 43, 47$.  But since $p=43$ and $p=47$ are congruent to 3
mod 4, the quadratic subfield of $\Q(\zeta_p)$ is imaginary for these
$p$.  On the other hand, since $2021>0$, the unique quadratic subfield
of $\Q(\RR_f)$ is real.  Thus we have linear disjointness for all $p$.
It follows that
$$\End(J^{(f_4,q)})=\Z[\zeta_q]$$
for all odd $p$.  
\end{sexs}

\begin{sthm}[Theorem 5.17 of \cite{ZarhinPisa}]
  Suppose that $m=3$, $p=2$, and $q=4$. Then $\End(J^{(f,q)})$ is an
  order in the matrix algebra of size two over the imaginary quadratic
  field $\Q(\zeta_4)=\Q(\sqrt{-1})$.
\end{sthm}

\begin{sthm}[Theorem 5.18 of \cite{ZarhinPisa}]\label{thm:cubic}
  Suppose that $k_0$ is an algebraically closed field of
  characteristic zero and $k=k_0(z)$ is the field of rational
  functions in the variable $z$. Suppose that $m=3$ or $4$ and
  $\Gal(f)=S_m$. If $(m,q) \ne (3,4)$ then
$$\End(J^{(f,q)})=\Z[\zeta_q].$$
\end{sthm}

\begin{sex}\label{ex:cubic}
Let $k_0=\bar{\Q}\subset \C$ and $f(x)=x^3-x-z$. One may view $\bar{\Q}(z)$ as a
subfield of $\C$, sending $z$ to any transcendental complex number $\beta$.
Then we get the complex polynomial $f_{3,\beta}(x)=x^3-x-\beta \in \C[x]$ and
conclude that for each $q \ne 4$ the complex abelian variety
$J^{(f_{3,\beta},q)}$ satisfies
$$\End(J^{(f_{3,\beta},q)})=\Z[\zeta_q].$$
\end{sex}

\section{Homomorphisms between superelliptic Jacobians}

We can apply the results of the previous section to obtain pairs of
superelliptic Jacobians with no homomorphisms between them.

\subsection{The case $g=y^n-a$}
In this section, we assume that $n\ge2$ and $g(y)=y^n-a$ for some non-zero
$a$ in $k$.  By Section~\ref{ss:3-monomials}, for all powers $q$ of a prime $p$
with $q>2$, $J_{g,q}$ is an abelian variety of CM type.

\begin{sthm}\label{thm:g-CM}
  Let $q=p^r$ for a prime number $p$.  Suppose that $m \ge 4$, $f$ is
  an irreducible polynomial over $k$ of degree $m$, and the 
  Galois group $\Gal(f)=S_m$ or $A_m$.  If $m=4$ then
  we also assume that $p$ is odd,  $\Gal(f)=S_4$, and at least one of the following two
  conditions holds:

\begin{itemize}
\item[(i)]  the unique quadratic subextension of
  $k(\RR_f)$ does not lie in $k(\zeta_p)$.

\item[(ii)] $k=\Q$ and $p$ is unramified in $\Q(\RR_f)$.
\end{itemize}

Then $\Hom(J(C_{f,q}), J(C_{g,q}))=\{0\}.$
\end{sthm}

\begin{proof}
  It suffices to check that $\Hom(J^{(f,p^i)},J(C_{g,q}))=\{0\}$ for
  all positive $i\le r$.  By Theorems~\ref{bign}, \ref{n4} and
  \ref{n4bis}, $J^{(f,p^i)}$ is absolutely simple and not of CM type
  over $\bar{k}$. But $J(C_{g,q})$ is of CM type over $\bar{k}$.
\end{proof}

\subsection{The case $m>n\ge4$}

\begin{sthm}\label{thm:m>n}
  Suppose that $m > n \ge 4$, the Galois groups $\Gal(f)=S_m$ or $A_m$
  and $\Gal(g)=S_n$ or $A_n$.  If $n=4$ then we also assume that $p$
  is odd, $\Gal(g)=S_4$, and at
  least one of the following two conditions holds:
\begin{itemize}
\item[(i)]  the unique quadratic subextension of $k(\RR_g)$ does not lie in $k(\zeta_p)$.
\item[(ii)] $k=\Q$ and $p$ is unramified in $\Q(\RR_g)$.
\end{itemize}
Then $\rk\Hom(J_{f,q},J_{g,q})$ is bounded by a constant independent
of $p$ and $q$.
\end{sthm}

\begin{proof}
It will suffice to show that 
$$\Hom(J^{(f,q_1)},J^{(g,q_2)})=\{0\}$$
for all divisors $q_1$ and $q_2$ of $q$ except possibly when
$q_1=q_2$ and $q_1|m$.  Under the hypotheses of the theorem, 
Theorems~\ref{bign}, \ref{n4} and \ref{n4bis} imply that $J^{(f,q_1)}$ and
$J^{(g,q_2)}$ are absolutely simple with endomorphism algebras
$\Q(\zeta_{q_1})$ and $\Q(\zeta_{q_2})$ respectively.  Thus either
$$\Hom(J^{(f,q_1)},J^{(g,q_2)})=\{0\}$$
or $J^{(f,q_1)}$ and $J^{(g,q_2)}$ are isogenous.  But if they were
isogenous, comparing endomorphism algebras we would have
$\phi(q_1)=\phi(q_2)$ and therefore $q_1=q_2$.  Comparing dimensions
(cf.~Section~\ref{ss:superelliptic-jacobians}), we would also have
$n=m-1$ and $q_1|m$.  Thus
$$\Hom(J^{(f,q_1)},J^{(g,q_2)})=\{0\}$$
unless $q_1=q_2$ and $q_1|m$.
\end{proof}

\section{Berger curves}\label{s:Berger}

\subsection{Notation}
Let $K=\kbar(t)$ be the
rational function field over $\kbar$.  Recall that given two rational
functions $f$ and $g$ on $\P^1$ over $k$, under mild hypotheses Lisa
Berger's construction \cite{Berger08} gives a smooth proper curve
$X_{f,g}$ over $K$ which is a model of the curve
$$\{f-tg=0\}\subset \P^1_K\times\P^1_K$$
Again under mild hypotheses, Berger computes the genus of $X_{f,g}$,
shows that $X_{f,g}$ is absolutely irreducible, and that in a suitable sense
it is associated to a tower of surfaces that are dominated by a
product of curves.  

More precisely, applying \cite{Berger08}*{Thm.~3.1} we have:

\begin{thm}
  Let $k$ be a field of characteristic zero and let $K=\kbar(t)$.  Let
  $f(x)$ and $g(y)$ be polynomials over $k$ of degrees $m$ and $n$ respectively with
  simple roots.  Let $X_{f,g}$ be the smooth proper model of
$$\{f-tg=0\}\subset \P^1_K\times\P^1_K$$
over $K$.  Then $X_{f,g}$ is absolutely irreducible and has genus
$$g_X=\frac{(m-1)(n-1)-\gcd(m,n)+1}{2}.$$
\end{thm}

Let $J$ be the Jacobian of $X_{f,g}$, an abelian variety of dimension $g_X$
over $K$.  For each positive integer $d$, let $K_d=\kbar(t^{1/d})$.
Building on \cite{Berger08}, Ulmer gives a formula in \cite{UlmerDPCT}
for the rank of the Mordell-Weil group $J(K_d)$ in terms of
homomorphisms between $J_{f,d}$ and $J_{g,d}$.

More precisely, we have:

\begin{thm}\label{prop:rank}
With notation as above, for all $d$ the $K_d/\kbar$-trace of $J$ is
zero.  Moreover, we have
$$\rk J(K_d)=\rk\Hom(J_{f,d},J_{g,d})^{\mu_d}-c_1d+c_2(d).$$
Here the exponent signifies those homomorphisms commuting with the
natural action of $\mu_d$ on both factors, $c_1$ is a non-negative integer, and
$c_2$ is a periodic function.  We have
$$c_2(d)=(m-1)(n-1)+\gcd(m,n,d)-1.$$
If for some $d$ strictly larger than $c_2(d)$ we have
$\Hom(J_{f,d},J_{g,d})^{\mu_d}=0$, then $c_1=0$.
\end{thm}

\begin{proof}
The vanishing of the $K_d/\kbar$ trace follows from
\cite{UlmerDPCT}*{5.6}.  The formula for the rank is \cite{UlmerDPCT}*{6.2} and
the definition of $c_2$ in \cite{UlmerDPCT}*{6.1.1} immediately gives the
expression above.  Since $c_1$ is non-negative and the right hand side
of the rank formula is also non-negative, if the $\Hom$ group is zero
for a large value of $d$, we must have $c_1=0$.
\end{proof}

\section{Bounded ranks}
We now assemble the pieces to give several examples of Jacobians over
$k(t)$ with bounded ranks in the tower $\kbar(t^{1/p^r})$.  We treat
only the most straightforward examples and the reader will have no
trouble seeing that there are several other ways to apply the basic
results of Sections~\ref{s:endos} and \ref{s:Berger}.

\subsection{Proof of Theorem~\ref{thm:main}}
For an integer $g_X\ge2$, let $m=2g_X+1$, $f(x)=f_m(x)=x^{m}-x-1$,
$n=2$, and $g(y)=y^2-1$.  Then Berger's curve $X_{f,g}$ is a smooth,
proper model of the hyperelliptic curve
$$ty^2=x^{2g_X+1}-x-1+t$$
of genus $g_X$. Let $J$ be the Jacobian of $X_{f,g}$.  Taking the
ground field $k$ to be $\Q$, we know (see Section~\ref{ss:endo-intro})
that $f$ has Galois group $S_m$.  Let $p$ be a prime number and $q$ a
power of $p$.  Applying Theorem~\ref{thm:g-CM}, we have that
$\Hom(J_{f,q},J_{g,q})=0$ and therefore by
Theorem~\ref{prop:rank}, the $\Qbar(t^{1/q})/\Qbar$-trace of $J$
is zero and 
$$\rk J(\Qbar(t^{1/q}))=(m-1)(n-1)=2g_X.$$

Note that $X_{f,g}$ is a quadratic twist of the hyperelliptic curve
$$X^{\prime}:\quad y^2=x^{2g_X+1}-x-1+t$$
and therefore $J$ is isomorphic over $\overline{\Q(t)}$ to the Jacobian
$J^{\prime}$ of $X^{\prime}$;  in particular, $\End(J^{\prime})=\End(J)$.
Clearly,  the Galois group of $x^{2g_X+1}-x-1+t$ over $\Q(t)$ coincides with
the Galois group of $x^{2g_X+1}-x-t$ over $\Q(t)$ and
therefore equals $S_m$ (even over $\bar{\Q}(t)$)
\cite{SerreGalois}*{Thm.~4.4.1 on p.~39}. The case $q=2$ and $k=\Q(t)$
of Theorem~\ref{bign} implies $\End(J^{\prime})=\Z$ and therefore $\End(J)=\Z$.
This shows in
particular that $J$ is absolutely simple.

It follows from Theorem~5.18(i) on p.~384 of \cite{ZarhinPisa} that in fact
the $\overline{\Q(t)}/\bar{\Q}$-trace of $J$ is zero.

This completes the proof of Theorem~\ref{thm:main}.  \qed

\subsection{}
Essentially identical arguments apply with $f(x)=x^{2g_X+2}-x-1$ and
$g=y^2-1$ and show that the Jacobian $J$ of the hyperelliptic curve 
$$ty^2=x^{2g_X+2}-x-1+t$$
is absolutely simple, has trivial $\Qbar(t^{1/q})/\Qbar$-trace, and has rank
$$\rk J(\Qbar(t^{1/q}))=2g_X+\gcd(q,2)-1.$$

The vanishing of the $\overline{\Q(t)}/\bar{\Q}$-trace of $J$ follows
from Theorem~1.2 of \cite{ZarhinSH} combined with Lemma~3.4 on p.~369
of \cite{ZarhinPisa}.

\subsection{A superelliptic generalization}
Consider the case $n\ge2$, $g=y^n-1$, $m\ge5$, $f=x^m-x-1$, and
$k=\Q$.  Berger's curve $X_{f,g}$ is the superelliptic curve
$$ty^n=x^m-x-1+t$$
of genus 
$$g_X=\frac{(m-1)(n-1)-\gcd(m,n)+1}{2}.$$ 
Let $J$ be the Jacobian of $X_{f,g}$ and let $q$ be a power of a prime
number $p$.  Applying Theorem~\ref{prop:rank} and
Theorem~\ref{thm:g-CM}, we have that
$$\rk J(\Qbar(t^{1/q}))=(m-1)(n-1)+\gcd(m,n,q)-1.$$
Note that this differs from $2g_X$ by an amount bounded independently of
$q$.

\subsection{More examples}
Now consider the case where $m$ and $n$
satisfy $m>n>4$ and $f$, $g$ and $p$ satisfy the hypotheses of
Theorem~\ref{thm:m>n}.  Berger's curve $X_{f,g}$ is as usual the
smooth proper model of $f-tg=0$.  It has genus
$$g_X=\frac{(m-1)(n-1)-\gcd(m,n)+1}{2}.$$ 
Let $J$ be the Jacobian of $X_{f,g}$.  For every
power $q$ of $p$, we have
$$\rk J(\kbar(t^{1/q}))\le(m-1)(n-1)+\gcd(m,n,d)-1+\epsilon$$
where $\epsilon$ is a constant which is independent of $p$ and $q$.

\subsection{Infinitely many non-isogenous examples over a Hilbertian
  field}
Let us discuss a generalization of previous examples that is based on
the notion of Morse polynomial \cite{SerreGalois}*{p.~39}. (Notice
that $x^m-x$ is a Morse polynomial and $f_m(x)=(x^m-x)-1$.)  Namely,
let $h(x)\in k[x]$ be a Morse polynomial of degree $m\ge 5$. Then the
Galois group of the polynomial of $h(x)-t$ over $k(t)$ is $S_m$
\cite{SerreGalois}*{Thm.~4.4.1 on p.~39}. Let us consider the subset
$A(h) \subset k$ consisting of all $c \in k$ such that the Galois
group of $h(x)-c$ over $k$ is $S_m$.   Assuming $c \in A(h)$, let
$f(x)=h(x)-c$ and $g(y)=y^2-1$.  The hyperelliptic curve 
$$X_{h-c, g}:\quad h(x)-c-t(y^2-1)=0$$
which is a $k(\sqrt{t})/k(t)$-twist of
$$X'_c:\qquad y^2 =h(x)-c+{t}=h(x)-(-t+c)$$
has genus $[\frac{n-1}{2}]$.  Clearly, the Galois group of
$h(x)-(-t+c)$ over $k(t)$ is also $S_m$; in particular, by
Theorem~\ref{bign}, the Jacobian $J_c$ of $X_{h-c,g}$ has zero
$\kbar(t)/\kbar$-trace and $\End(J)=\Z$.  Arguing as in the proof of
Theorem~\ref{thm:main}, we have that the rank of
$J_c(\kbar(t^{1/p^r}))$ is bounded independently of $p$ and $r$.  Let
$J'_c$ be the Jacobian of $X'_c$ and note that $J'_c$ and $J_c$ become
isomorphic over $k(t^{1/2})$.

Now suppose that $d$ is a different element of $A(h)$; assume in
addition that the splitting field $k_c$ of $h(x)-c$ and the splitting
field $k_d $ of $h(x)-d$ are linearly disjoint over $k$.  Let us
consider the Jacobians $J'_c$ and $J'_d$ as abelian varieties over the
complete discrete valuation field $k((t))$.  Clearly, they both have
good reduction at $t=0$; their reductions are the Jacobians of
$y^2=h(x)-c$ and $y^2=h(x)-d$ respectively. In addition, the field of
definition of points of order $2$ on $J'_c$ (resp.~on $J'_d$)
coincides with $k_c((t))$ (resp.~with $k_d((t))$; in particular, those
splitting fields are linearly disjoint over $k((t))$. Now it follows
from Theorem~1.2 of \cite{ZarhinSH} that $J'_c$ and $J'_d$ are
\emph{not} isogenous over $\overline{k((t))}$!  Therefore $J_c$ and
$J_d$ are not isogenous over $\overline{k(t)}$.

Notice that if $k$ is Hilbertian (e.g., a number field) then Hilbert's
irreducibility theorem guarantees the existence of plenty of such $c$
and $d$.

\subsection{One more elliptic example}
As in Example~\ref{ex:cubic}, let $k=\C$ and let $\beta$ be a
transcendental complex number.  If we put $g(x)=f(x)=x^3-x-\beta$ then
Berger's construction gives us an an elliptic curve
$X_{f,g}:(x^3-x-\beta) - t(y^3-y-\beta)$ over the function field
$\C(t)$ that appears (in slightly disguised form) in
\cite{UlmerDPCT}*{Sect.~9}.  Applying Theorem~\ref{thm:cubic} and the
rank formula in \cite{UlmerDPCT}*{Section~9.3} we have that the rank
of $X_{f,g}(\C(t^{1/p^r}))$ is bounded for all primes $p$ and all
$r\ge1$.

\section*{Acknowledgements}
Ulmer's research was partially supported by NSF grant DMS
  0701053.

\end{document}